\documentclass[11pt]{smfart}
\usepackage{verbatim,amssymb,url,xspace,smfenum}
\RequirePackage[T1]{fontenc}
\RequirePackage{lmodern}
\RequirePackage[francais,english]{babel}
\usepackage[latin1]{inputenc}

\let\ba\smallsetminus
\let\le\leqslant
\let\ge\geqslant
\let\leq\leqslant
\let\geq\geqslant
\let\hat\widehat
\let\dpl\displaystyle
\let\wh\widehat

\let\ov\overline

\let\ldots\dots
\def\sfrac#1#2{{#1}/{#2}}
\def\resp{resp.\kern.3em}
\newcommand{\ptbl}{.\kern .15em }
\def\bbC{\mathbb{C}}
\def\bbR{\mathbb{R}}
\def\bbZ{\mathbb{Z}}
\def\bbT{\mathbb{T}}
\def\cF{\mathcal{F}}
\def\cN{\mathcal{N}}
\newcommand{\reel}{\mathop{\mathrm{R\acute{e}}}\nolimits}

\theoremstyle{plain}
\newtheorem*{theoreme*}{\theoname}

\newenvironment{enumeratei}
{\bgroup\def\theenumi{\roman{enumi}}\def\theenumii{\arabic{enumii}}\begin{enumerate}}
{\end{enumerate}\egroup}
\newenvironment{enumeratea}
{\bgroup\def\theenumi{\alph{enumi}}\begin{enumerate}}
{\end{enumerate}\egroup}

\begin{document}
\frontmatter
\title{Analyse et synthèse harmoniques}
\author[J.-P. Kahane]{Jean-Pierre Kahane}
\address{Laboratoire de Mathématique d'Orsay, UMR 8628\\
Bâtiment 425\\
Faculté des Sciences d'Orsay\\
Université Paris-Sud 11\\
F-91405 Orsay Cedex}
\email{Jean-Pierre.Kahane@math.u-psud.fr}

\thanks{Texte des exposés faits aux Journées mathématiques X-UPS le 10 mai 2011}
\maketitle
\mainmatter

\part*{Coup d'{\oe}il sur l'analyse de Fourier}

Embrasser l'analyse de Fourier d'un coup d'{\oe}il est hors de ma portée. Je promènerai le regard sur l'histoire et sur certains des termes en usage. L'histoire est ancienne et elle est instructive. Je commencerai par Platon.

Chez Platon, les mathématiques comprennent cinq parties: les nombres, les figures planes, les solides, la musique et l'astronomie. La musique et l'astronomie sont comme deux s{\oe}urs, liées par la musique des sphères.
Et le programme astronomique de Platon est de rendre compte du mouvement des astres errants, les planètes, pour le rendre conforme à l'harmonie du monde et à la dignité des dieux.

Le système de Ptolémée a réalisé ce programme, par une superposition de mouvements circulaires uniformes, le grand cycle et les épicycles. Et quoique le système de Ptolémée soit abandonné, la décomposition d'un mouvement en mouvements périodiques et l'utilisation des fonctions trigonométriques a été une constante de l'astronomie jusqu'à nos jours. Aujourd'hui l'astrophysique double et dépasse l'astronomie dans l'approche des phénomènes cosmiques. Mais la perle de l'astronomie qu'est la découverte des exoplanètes repose sur l'analyse de Fourier des spectres émis par les étoiles. 

Quant à la musique, son lien aux mathématiques est multiple. Le principal est la décomposition d'un son en harmoniques et sa reconstitution comme somme d'harmoniques. C'est au sens propre l'analyse et la synthèse harmoniques. A cela s'attachent les gammes, la distinction des timbres, les cordes vibrantes et la grande controverse qui au 18\ieme siècle a opposé Daniel Bernoulli à d'Alembert, Euler et Lagrange sur la possibilité de représenter une fonction par une série trigonométrique. Aujourd'hui l'analyse par ondelettes complète ou supplante l'usage des fonctions trigonométriques pour tenir compte des échelles de temps aussi bien que des fréquences, et le champ de l'analyse de Fourier s'étend par ses méthodes comme par ses usages.

Fourier, selon Riemann, est le premier à avoir compris complètement la nature des séries trigonométriques, en associant ce que j'appellerai l'analyse, les formules intégrales donnant les coefficients, et la synthèse, la représentation d'une fonction par une série trigonométrique. C'est à cause de Riemann que nous parlons aujourd'hui de séries de Fourier. En France, Fourier a été longtemps méconnu. Arago, dans son éloge de Fourier, vante grandement le savant et le politique, explique l'importance de sa théorie analytique de la chaleur, mais ne dit pas un mot des séries trigonométriques, c'est-à-dire de l'outil que Fourier a forgé pour calculer des solutions des équations intégrales de la chaleur. Fourier y accordait une grande importance et, contrairement à une idée reçue, il a développé cet outil en véritable théorie, exemples, applications et démonstrations, avec une grande rigueur. Mais il s'est heurté à l'incompréhension persistante de Lagrange, et deux pages dans les manuscrits de Lagrange, que j'ai consultées et commentées, confirment que Fourier avait raison contre Lagrange. Mais Lagrange était à l'époque de Fourier le plus respecté des mathématiciens français, et son jugement négatif sur Fourier a traversé les siècles. Dans l'édition que je possède d'Encyclop{\ae}dia Universalis il n'y a pas d'article sur Fourier. 

Fourier attachait une portée universelle à ses formules. Il précisait bien, et le premier, que la donnée d'une fonction était celle de son domaine de définition en même temps que d'une loi ou une figure ; et, sans utiliser ces termes, il distinguait soigneusement les intervalles ouverts et les intervalles fermés. Mais, après avoir multiplié les exemples et donné quelques preuves, il s'était aventuré à dire que toute fonction était sujette à cette analyse et représentable par une série trigonométrique qui converge vers la fonction. Littéralement c'est faux. 
Pour appliquer les formules intégrales, il faut que la fonction soit intégrable ; et la convergence des séries est un sujet difficile. Les premiers pas pour éclaircir la question sont dus à Dirichlet, avec le premier théorème général de convergence et le premier exemple de fonction non intégrable dans la conception de l'époque. Tous les concepts d'intégration, à commencer par l'intégrale de Riemann, sont liés aux formules de Fourier et à leurs conditions de validité. Au delà de Riemann, on pense à Lebesgue, Denjoy, Laurent Schwartz.
Quant à la convergence, elle est mise en question par les contre-exemples: une fonction continue dont la série de Fourier diverge en un point (du Bois-Reymond) ou même sur un ensemble donné de mesure nulle (Kahane-Katznelson), et une fonction intégrable au sens de Lebesgue dont les série de Fourier diverge partout (Kolmogorov). Le théorème de convergence de Carleson, inattendu à son époque (1966), dit que la convergence a lieu presque partout quand la fonction est de carré intégrable ; on peut améliorer cette condition, sans parvenir bien sûr aux fonctions intégrables.

De même que le concept d'intégrale, c'est le concept de série qui est en cause. Pourquoi s'attacher à la convergence, qui a d'ailleurs plusieurs significations dès qu'on passe à plusieurs variables et à des séries multiples, alors qu'il y a des procédés de sommation utilisables ? Le tournant est pris avec le théorème de Féjer de 1900: pour les fonctions continues, les moyennes arithmétiques des sommes partielles convergent. Elles convergent même uniformément. La convergence dans les espaces fonctionnels apparaît peu après ; les fonctions sont des points, les sommes partielles d'autres points, et l'espace des fonctions de carré intégrable s'impose à l'attention avec la formule de Parseval, établie dans ce cadre par Fatou, et surtout le théorème de Riesz-Fischer, qui établit l'isomorphisme isométrique de $L^2$ et de $\ell^2$ par transformation de Fourier (1907).

Fischer et Frédéric Riesz, indépendamment, ont pensé à une géométrisation des espaces de fonctions, et le lemme fondamental pour la preuve de leur théorème est le même ; nous l'exprimons aujourd'hui en disant que $L^2$ est complet. Mais cette formulation ne date que de Banach dans sa théorie des opérations linéaires de 1930. Il fallait jusque là une longue phrase pour le dire. C'est un exemple où les définitions, tardives, viennent exprimer le suc d'une méthode, avant de servir de base à de nouveaux développements.

Autre exemple, toujours tiré de l'analyse de Fourier. Au lieu de $L^2$ et $\ell^2$, Norbert Wiener s'est attaché à $L^1$ et $\ell^1$. Leurs transformées de Fourier sont un champ d'étude toujours ouvert, avec des applications surprenantes en théorie du signal (Donoho, Candès etc.). Les problèmes d'analyse et de synthèse y prennent un aspect différent, plus algébrique: dans l'algèbre des fonctions sommes de séries trigonométriques absolument convergentes, les fonctions nulles sur un ensemble donné forment un idéal fermé ; y en a t-il d'autres ?
C'est bien la cas, comme Malliavin l'a montré en 1959. Le point de départ est l'algèbre de Wiener, qui se voit soit comme l'algèbre multiplicative des fonctions sommes de séries trigonométriques absolument convergentes, soit comme algèbre de convolution $\ell^1$. La dernière expression est plus rapide, et la notion de convolution, si fondamentale en analyse, la rend très parlante. Mais en 1930, et même quand Laurent Schwartz a élaboré sa théorie des distributions, on ne parlait pas encore de convolution ; c'était Faltung en allemand, produit de composition en français. Dans le traité de Widder \og Laplace transforms\fg, qui date de 1941 , c'est comme \og Stieltjes transforms\fg que sont introduites les convolutions de mesures, avec une note en bas de page qui indique le terme de convolution comme une timide nouveauté parallèlement au terme de Faltung, aussi utilisé en anglais. La convolution apparaît partout, mais c'est Wiener qui a dégagé son caractère fondamental , et c'est pourquoi Faltung est devenu pour un temps l'expression de la notion. Aujourd'hui il est raisonnable de placer la convolution au départ d'un cours d'analyse de Fourier.

S'agissant du vocabulaire, comment situer l'analyse harmonique ? C'est un champ largement ouvert sur l'ensemble des mathématiques. Historiquement, on vient de voir ses relations avec les équations différentielles de la mécanique céleste, avec les équations aux dérivées partielles des cordes vibrantes et de la chaleur, avec la théorie des fonctions d'une variable réelle, et on a évoqué ses relations actuelles avec la statistique et le traitement des données. Le lien aux probabilités est ancien et profond. Le cadre de l'analyse harmonique commutative est celui des groupes abéliens localement compacts, et la dualité entre ces groupes est exprimée par la transformation de Fourier. La thèse de Tate a montré l'importance de cette approche en théorie des nombres. L'analyse harmonique non commutative est liée aux groupes non commutatifs et à leurs représentations. L'analyse harmonique abstraite part de la théorie de Gelfand des anneaux normés, autrement dits algèbres de Banach. L'analyse de Fourier classique traite de transformations intégrales, transformation de Fourier d'abord, et aussi intégrales singulières. Ses usages dans toutes les sciences ont été multipliés par la transformation de Fourier rapide puis par les ondelettes et leurs variantes. Il est clair que toute définition de l'analyse harmonique en serait une limitation injustifiée. Mais outre son étendue, on peut repérer des questions, des méthodes, des théories qui, elles, peuvent être identifiées et formalisées.

Les formules de Fourier sont l'exemple de base. Leur énoncé par Fourier exprime leur généralité, de façon formellement incorrecte. Elles ne constituent pas un théorème, mais elles ont engendré des théorèmes et permis de définir d'importantes notions. Mieux qu'un théorème, elles ont constitué un programme, à savoir, donner des conditions de leur validité. Elles ont ensuite constitué un paradigme pour tous les développements orthogonaux. Une raison de leur succès est sans doute qu'elles établissent un pont entre deux classes d'objets, en l'occurrence des fonctions et des suites. La dualité de Fourier est un modèle de traduction d'un langage dans un autre, un trait commun à d'autres grands programmes.

\enlargethispage{2\baselineskip}%
Le coup d'{\oe}il pourrait se poursuivre, aussi bien sur l'histoire que sur l'actualité de l'analyse de Fourier. Ce sera pour une part l'objet des exposés suivants. En attendant, voici quelques références.

René Spector (1970), dans l'article \og Harmonique (analyse)\fg d'\emph{Encyclop\ae dia Universalis}, donne un bon aperçu en quatre parties: séries de Fourier, analyse et synthèse harmonique, transformation de Fourier, groupes commutatifs localement compacts, avec un mot sur les groupes de Lie.

Jean-Paul Pier (1990), dans \emph{L'analyse harmonique, son développement historique} (Masson) prend pour point de départ l'apparition de la théorie des groupes topologiques généraux, de la mesure de Haar et des algèbres de Banach, donc il se limite pour l'essentiel au \textsc{xx}\ieme siècle, avec beaucoup de détails sur les questions qu'il traite, et une abondante bibliographie (40 pages) \cite{Pier90}.

Elias Stein et Rami Shakarchi, dans \emph{Fourier analysis, an introduction} (2002) ont un point de vue très différent. Au départ, pas d'intégrale de Lebesgue ni de groupe: intégrale de Riemann, droite et cercle suffisent, puis $\bbR^d$, en enfin~$\bbZ_N$, avec beaucoup d'applications \cite{S-S03}. C'est le premier livre d'une série qui se poursuit avec l'analyse complexe, la théorie de la mesure, l'intégrale de Lebesgue et les espaces de Hilbert, puis les distributions de Schwartz et les probabilités, dans un foisonnement d'applications.

On voit la variété des approches. D'excellents livres ont paru dans les années 1930. Je signale seulement le grand classique, le traité d'Antoni Zygmund \emph{Trigonometrical series} qui, dans la seconde édition (1959) est devenu \emph{Trigonometric series}; la \og troisième\fg édition (2002) reproduit la seconde \cite{Zygmund1, Zygmund2, Zygmund3, Zygmund4}.

\emph{Abstract harmonic analysis} de L\ptbl Loomis (1953) a permis aux gens de ma génération de s'initier aux algèbres de Banach (les anneaux normés d'Israel Gelfand). Le seul point commun aux livres de Zygmund et de Loomis est la théorie de Norbert Wiener des transformées de Fourier des fonctions intégrables \cite{Loomis53}.

Sur ce sujet, qui a fait de grands progrès au cours des années 1950, Walter Rudin fait le point avec \emph{Fourier transforms on groups} (1962) \cite{Rudin1,Rudin2}. Des progrès ultérieurs sont exposés dans \emph{Essays in commutative harmonic analysis} (1980) de Colin Graham et O.C\ptbl McGehee \cite{G-Mc79}.

Henry Helson (1983) titre simplement son livre \emph{Harmonic analysis} et il se centre sur un petit nombre de questions, typiques et essentielles; son but, dit-il, est la simplicité, et c'est un modèle d'écriture \cite{Helson1, Helson2}.

Au regard de H\ptbl Helson, Tom Körner est foisonnant dans ses deux livres \emph{Fourier Analysis} (1988) et \emph{Exercises} (1993). Ces deux livres, savants et riches, pleins d'aperçus et d'anecdotes, offrent une délectable promenade dans l'analyse de Fourier \cite{Korner88, Korner93}.

C'est aussi une promenade historique que Pierre-Gilles Lemarié-Rieusset et moi proposons dans \emph{Séries de Fourier et ondelettes}, en citant abondamment les grands ancêtres, et en détaillant l'histoire récente des ondelettes \cite{KaLe}.

Je termine par le livre d'Yitzhak Katznelson, \emph{An introduction to harmonic analysis}, parce que j'ai assisté à sa naissance dans le Languedoc et à son audience croissante aux États-Unis et dans le monde (1968--1976--2006). C'est beaucoup plus qu'une introduction, c'est l'exposé clair et vivant de l'analyse harmonique par l'un de ses principaux acteurs \cite{Katznelson1, Katznelson2, Katznelson3}.

\part{Statistique et Fourier, dans~l'histoire~et~aujourd'hui}

Ce sera une promenade historique, dans laquelle je pense m'attarder sur Fourier, mais qui aura comme point de départ la méthode des moindres carrés et comme point d'arrivée l'échantillonnage parcimonieux, alias \emph{compressed sensing}.

\bigskip
\centerline{\textbf{Personnages}}

\begin{multicols}{3}
\noindent Clairaut 1713-1765\\
Lagrange 1736-1813\\
Laplace 1749-1827\\
Legendre 1736-1813\\
Fourier 1768-1830\\
Gauss 1777-1855\\
Riemann 1826-1866\\
E. Fischer 1875-1854\\
Lebesgue 1876-1941\\
F. Riesz 1880-1948\\
N. Wiener 1894-1964\\ \\
Yves Meyer\\
Stephane Mallat\\
David Donoho\\
Emmanuel Candès\\
JPK
\end{multicols}

Le personnage principal sera Fourier, mais je parlerai aussi de ses aînés, Legendre et Laplace, et de ses successeurs, parmi lesquels Riemann et Lebesgue. Dans la foison des théories et des résultats du 20\ieme siècle j'insisterai sur deux aspects: ce qui concerne~$\ell^2$, l'espace des suites de carré sommable, qui a été mis en évidence par Frédéric Riesz et Ernst Fischer, et ce qui concerne~$\ell^1$, l'espace des suites sommables, qui a été découvert par Norbert Wiener. La statistique sera l'occasion de voir~$\ell^2$ et~$\ell^1$ s'introduire naturellement, avec les moindres carrés de Legendre et compagnie, et le \emph{compressed sensing} d'Emmanuel Candès, David Donoho et compagnie.
\enlargethispage{\baselineskip}%

Pour ne pas risquer de manquer de temps en m'attardant en route, je commence par la fin du voyage: le théorème de Candès, Romberg et Tao de 2006 sur la reconstitution d'un signal ou d'une image à partir de données apparemment insuffisantes. En théorie du signal, on enseigne que pour reconstituer une fonction dont la transformée de Fourier est portée par une bande de largeur~$L$ on~a besoin d'un échantillonnage sur une progression~arithmétique de raison~$a/L$, $a$~dépendant de la normalisation. En limitant la fonction à un intervalle de longueur $H$, on~a donc besoin de $HL /a$ points. Or voici le théorème:

On sait qu'un signal est constitué de $T$ pics, sans connaître leur position, et on~a accès à sa transformée de Fourier. Il s'agit de reconstituer le signal à partir de la transformée de Fourier de la manière la plus économique. On identifie le signal à une fonction $f$ portée par $T$ points du groupe cyclique à~$N$ points, $N$ choisi assez grand. Ce groupe cyclique $\bbZ_N$ que je note~$G$, est isomorphe à son dual $\hat G$, le groupe des fréquences, sur lequel est définie la transformée de Fourier $\hat f$. On choisit au hasard, avec la probabilité naturelle, $\Omega$ fréquences $\omega$ et on considère les $\hat f(\omega)$. Supposons
\[
T < C \Omega / \log N\qquad\text{avec}\quad C = 1 / 23 (M+1).
\]
Alors on peut reconstituer $f$ à partir des $\hat f(\omega)$ avec une probabilité $1 - c / N^M)$, où $c$ est une constante absolue. Et la reconstitution se fait par un procédé d'analyse convexe, à savoir en minimisant la norme~$\ell^1$ des fonctions définies sur~$G$ dont la transformée de Fourier coïncide avec $\hat f$ sur les $\Omega$ fréquences choisies. On voit l'économie considérable: l'utilisation classique de la transformée de Fourier nécessite $N$ points, et l'échantillonnage aléatoire le réduit à un multiple de $T \log N$.

Rappelons que la norme~$\ell^1$ d'une fonction $h$ définie sur~$G$ est la somme des valeurs absolues des $h(g)$, $g$ parcourant~$G$. C'est bien une fonction convexe de $h$. Il y a aujourd'hui des procédés variés pour traiter de tels problèmes de minimisation, et je dirai un mot de leur histoire.

\medskip
Remontons à 200 ans en arrière, plus précisément en 1805. Legendre a écrit un mémoire intitulé \og Nouvelles méthodes pour la détermination des orbites des comètes\fg où il utilise un procédé qu'il décrit soigneusement dans un appendice de quatre pages intitulé \og Sur la méthode des moindres quarrés\fg, et il applique ensuite le même procédé à la mesure des degrés du méridien terrestre, un sujet de grand intérêt à l'époque.\enlargethispage{\baselineskip}%

Cet appendice est lumineux. En voici le début:

\og Dans la plupart des questions où il s'agit de tirer des mesures données par l'observateur les résultats les plus exacts qu'elles peuvent offrir, on est presque toujours conduit à un système d'équations de la forme
\[
E = a + bx + cy + fz + \text{etc.}
\]
dans lesquelles $a, b, c, f$, etc.\ sont des coefficients connus, qui varient d'une équation à l'autre, et $x, y, z$, etc.\ sont des inconnues qu'il faut déterminer par la condition que la valeur de $E$ se réduise, pour chaque équation, à une quantité ou nulle ou très petite.\fg

Je résume la suite. En général, il y a plus d'équations que d'inconnues, et il n'est pas possible d'annuler toutes les erreurs $E, E', E''$ etc. Le procédé consiste à rendre minimum la somme des carrés des erreurs. En dérivant la somme des carrés des erreurs par rapport à chacune des variables, on obtient des formes affines en nombre égal à celui des variables, et en les annulant on peut calculer ces variables, et aussi les erreurs $E, E', E''$ etc. Legendre précise qu'\og on~aura soin d'abréger tous les calculs, tant des multiplications que de la résolution, en n'admettant dans chaque opération que le nombre de chiffres entiers ou décimaux que peut exiger le degré d'approximation dont la question est susceptible.\fg

Si parmi les erreurs ainsi calculées il en est d'aberrantes, on rejette les équations correspondantes et on reprend le calcul, dont la majeure partie (les multiplications) aura déjà été faite.

Au prix de la multiplication des coefficients et de la résolution d'un système d'équations linéaires on~a donc un procédé général et efficace pour réaliser le programme indiqué par Legendre. Le procédé prend une forme remarquable quand il n'y a qu'une variable~$x$, qui représente un point sur la droite ou dans le plan ou l'espace, et des observations $x', x'', x'''$ etc. Le minimum de la somme des carrés distances au point~$x$ des points $x', x'', x'''$ etc.\ est atteint au contre de gravité de $x', x'', x'''$ etc. Comme le dit Legendre en conclusion,

\smallskip
\noindent
\og on voit donc que la méthode des moindres carrés fait connaître, en quelque sorte, le centre autour duquel viennent se ranger tous les résultats fournis par l'expérience, de manière à s'en écarter le moins qu'il est possible.\fg

\smallskip
Il y eut une querelle de priorité entre Legendre et Gauss sur la méthode des moindres carrés. Gauss l'a publiée plus tard et sous certaines hypothèses sur la loi des erreurs (en gros, qu'elles soient gaussiennes) et il l'a utilisée sans la publier avant 1805. Mais il semble bien que Legendre, et d'autres, l'aient aussi utilisée sans la formaliser avant 1805. Stigler, dans son histoire de la statistique, consacre son premier chapitre à la méthode des moindres carrés exposée par Legendre (en donnant en illustration l'image d'un autre Legendre, cela est une tout autre histoire), avec un commentaire enthousiaste:

\smallskip
\noindent
\og For stark clarity of exposition the presentation is unsurpassed; it must be counted as one of the clearest and most elegant introductions to a new statistical method in the history of statistics.\fg

\smallskip
Et surtout Stigler situe la méthode parmi celles qui avaient cours
précédemment, à savoir les combinaisons ad hoc des équations pour se ramener à autant d'équations que d'inconnues. Il peut valoir la peine d'examiner ces solutions \emph{ad hoc}. Pour autant que j'en puisse juger, elles combinent les seconds membres avec des coefficients qui sont toujours $1$, $0$ ou $-1$. L'avantage est d'éviter les multiplications. L'inconvénient est de faire des choix ad hoc au vu des équations. Maintenant qu'il est d'usage de faire des choix aléatoires, on peut se demander si des choix aléatoires de coefficients $1$, $0$ et $-1$, suivant une loi ne dépendant que du nombre d'équations et du nombre d'inconnues, ou plus simplement des choix aléatoires de $1$ et $-1$ avec la probabilité naturelle, ne donneraient pas des solutions presque aussi bonnes avec, par répétition du procédé, une estimation~assez correcte de l'erreur possible. Même s'il s'agit du centre de gravité, le choix aléatoire des coefficients $1$ et $-1$ à somme nulle donne une estimation rapide de l'erreur possible, au lieu du calcul de la variance.

Ce qui manquait à l'exposé de Legendre était justement l'estimation des erreurs possibles sur les quantités calculées. La liaison~avec les probabilités allait être faite par Gauss, et reprise par Laplace dès 1810 dans son \og mémoire sur les intégrales définies et leur application~aux probabilités\fg ; les applications comprennent la théorie des erreurs, celle de la combinaison des observations (dont je viens de parler), et la méthode des moindres carrés, dont Laplace déclare:

\smallskip
\noindent
\og Cette méthode, qui jusqu'à présent ne présentait que l'avantage de fournir sans aucun tâtonnement les équations finales nécessaires pour corriger les éléments, donne en même temps les corrections les plus précises.\fg

\smallskip
Bien avant Gauss, Laplace avait en main tous les éléments pour lier les probabilités et la méthode des moindres carrés. L'apparition de la fonction de Gauss et le calcul de son intégrale sur la droite apparaît chez lui dès 1774, dans son \og mémoire sur la probabilité des causes par les évènements\fg. En effet, la loi d'une grandeur x quand on connaît la loi des erreurs possibles et les résultats d'observation $x', x'', x'''$ est bien un calcul de probabilités des causes quand on connaît les évènements. Laplace en donne la forme générale
\[
\varphi(x'-x) \varphi(x''-x) \varphi(x'''-x)
\]
quand $\varphi (\cdot)$ est la loi des erreurs. Le cas particulier où $\varphi$ est la fonction de Gauss (Gauss est né en 1777) est aussi un cas limite, et la méthode des moindres carrés avec expression des corrections aurait pu apparaître dès ce moment, mais ce n'a pas été le cas. Laplace attribuait une grande valeur à la théorie des probabilités, il avait contribué au calcul de grandeurs à partir d'observations en combinant des équations, mais le lien entre les deux lui a échappé. Il a été découvert par Gauss avant que Laplace l'expose en 1810 puis dans son grand traité \og Théorie analytique des probabilités\fg.

La méthode des moindres carrés avec estimation des corrections était-elle pratiquée avant d'être parfaitement formalisée ? On peut se poser la question en lisant l'exposé très pédagogique qu'en fait Fourier dans deux mémoires tardifs, de 1826 et 1828, dans le second desquels il mentionne en passant l'avoir appliquée lorsqu'il était en Egypte, donc avant 1800, pour la détermination de la hauteur de la pyramide de Kéops. Ces mémoires ne contiennent pas de démonstrations, mais des énoncés et des applications, et ils sont destinés à des statisticiens. Depuis 1815, après avoir été démis de ses fonctions de préfet de l'Isère par Napoléon~au retour de l'ile d'Elbe, il s'était installé à Paris et avait été chargé, grâce au comte de Chabrol son~ancien élève, de la direction des services de statistique du département de la Seine.

Puis, après son élection tardive à l'Académie des sciences en 1817, il faisait figure de connaisseur en statistique. Son influence peut se mesurer par l'audience du prix Montyon de statistique dont il avait été le rapporteur pour l'adoption comme prix de l'académie, et plus encore par le témoignage de Quételet qui lui doit, dit-il, ses connaissances en probabilités.

C'est le moment de parler de la transformation de Fourier. A la suite des séries de Fourier, elle apparaît clairement dans l'ouvrage principal de Fourier, la Théorie analytique de la chaleur, qui sous sa forme achevée a été publiée en 1822. Elle était déjà là dans le mémoire de Fourier de 1811 qui lui a valu le Prix de l'Académie des sciences (ou plutôt, les académies royales ayant été supprimées par la Révolution, de la première classe de l'Institut de France) et qui n'avait pas été publié à l'époque. On y trouve la formule d'inversion et de nombreuses applications. Mais là encore, il y a des antécédents. Le plus important est la présentation par Clairaut de la transformation de Fourier discrète, celle dont il est question~au début de mon exposé, dans son mémoire de 1754 sur l'orbite apparente du soleil ; naturellement, ni Clairaut ni personne ne pouvait soupçonner alors l'importance de la notion, et l'influence de Clairaut sur Fourier et ses successeurs semble avoir été nulle.

Aujourd'hui le nom de Fourier est utilisé principalement pour ce qui faisait le plus question à son époque, les séries de Fourier, les intégrales de Fourier, la transformation de Fourier, et plus généralement l'analyse de Fourier. C'est à Riemann que l'on doit l'appellation de séries de Fourier pour les séries trigonométriques dont les coefficients sont donnés par les formules intégrales connues comme formules de Fourier. Dans sa thèse sur les fonctions représentables comme sommes de séries trigonométriques Riemann fait une belle étude historique du sujet, il en montre les antécédents, puis met en évidence le couplage entre les formules intégrales donnant les coefficients et l'expression de la fonction comme somme de série trigonométrique, en déclarant que Fourier est le premier à avoir compris la portée de ce couplage. Aujourd'hui nous pouvons exprimer ce couplage comme celui de l'analyse et de la synthèse. Fourier avait insisté sur sa valeur universelle, et là un tempérament est nécessaire: la théorie ne peut s'appliquer qu'à des fonctions intégrables, et alors la convergence de la série pose problème.

Ainsi on doit considérer les formules de Fourier, celles qui expriment les coefficients à partir de la fonction et celle qui exprime la fonction à partir des coefficients, comme un programme plutôt qu'un énoncé. Comme énoncé, il a amené à la recherche de contre-exemples instructifs. Mais comme programme, il est à l'origine de la formalisation de la notion d'intégrale, avec l'intégrale de Riemann, puis celle de Lebesgue, puis celle de Denjoy, puis des distributions de Laurent Schwartz. Il est aussi à l'origine des différents moyens d'exprimer la représentation de la fonction, ou du signal, comme somme d'une série trigonométrique: les procédés de sommation~au lieu de la seule considération des sommes partielles, et les espaces fonctionnels dans lesquels les fonctions apparaissent comme des points.

Parmi ces espaces fonctionnels le premier et le plus important est l'espace de Hilbert. On peut à volonté le considérer comme un~$\ell^2$, espace des suites de carré sommable, ou comme un~$L^2$, espace des fonctions de carré intégrable. On doit préciser quels sont les indices pour les suites, et de quelles fonctions il s'agit. Le cas classique, qui fait l'objet du théorème de Riesz-Fischer, est relatif à~$\ell^2(\bbZ)$ et~$L^2 (\bbT)$, $\bbT$ étant le cercle $\bbR/\bbZ$. L'isomorphisme est assuré par la transformation de Fourier: c'est la correspondance entre une fonction périodique (définie sur $\bbT$) et la suite de ses coefficients de Fourier. Voici des points de repère chronologiques. L'intégrale de Lebesgue date de 1901. L'égalité de Parseval, ainsi nommée par Fatou, se trouve dans sa thèse de 1906: elle exprime que la transformation de Fourier réalise une isométrie de~$L^2$ vers~$\ell^2$, mais non que cette isométrie est surjective. En 1907, Frédéric Riesz et Ernst Fischer, indépendamment ou plutôt concurremment, ont publié deux séries de notes aux Comptes rendus de l'Académie des sciences établissant le théorème et jetant les bases de ce que Fischer a appelé une géométrie des fonctions. La priorité n'est pas claire, et on peut aussi bien parler du théorème de Fischer-Riesz. La clé pour passer de Parseval à Riesz-Fischer est le fait que~$L^2$ est un espace complet. Et ce dernier fait montre la supériorité de l'intégrale de Lebesgue sur celle de Riemann.

Cependant raconter l'histoire ainsi est un renversement de son cours réel. En 1907, on ne disposait pas des symboles~$L^2$ ni~$\ell^2$. La notion d'espace normé n'existait pas encore, et celle d'espace normé complet encore moins. Il fallait de longues périphrases pour parler de suites de Cauchy, et même pour énoncer le théorème. C'est un cas typique ou l'outil de la démonstration, le fait que~$L^2$ est complet, est devenu le théorème principal, et que pour l'énoncer en trois mots on~a eu besoin de définir les termes aujourd'hui classiques de~$L^2$ et d'espace complet: l'énoncé sous cette forme n'apparaît qu'en 1930 dans la Théorie des opérations linéaires de Banach. En bref, le lemme est devenu théorème et sa substance est passée dans les définitions.

Le cadre général est apparu au cours des années 1930. C'est celui des groupes abéliens localement compacts et de leur dualité. Ainsi $L^2 (\bbT)$ et~$\ell^2(\bbZ)$ sont en dualité, de même que~$L^2(\bbR)$ avec lui-même ou~$\ell^2(\bbZ/ N\bbZ)$ avec lui même. La transformée de Fourier réalise un isomorphisme isométrique entre~$L^2 (G)$ et~$L^2 (\hat G)$,~$G$ et $\hat G$ étant deux groupes abéliens localement compacts en dualité. La notation~$\ell$ au lieu de~$L$ est utilisée lorsque le groupe est discret.

L'aspect géométrique est traduit aussi par la terminologie des bases orthonormales et des systèmes orthogonaux. Avec leur ancêtre, la méthode des moindres carrés, et tous les développements d'analyse liés aux systèmes orthogonaux, les suites de carré sommable ont une place de choix dans l'univers mathématique.

Et cependant il y a plus naturel encore que~$L^2$ et~$\ell^2$. C'est~$L^1$ et~$\ell^1$. En termes classiques, les transformées de Fourier fournissent les séries de Fourier-Lebesgue (les coefficients sont donnés par les formules de Fourier avec comme intégrale l'intégrale de Lebesgue), appelées aussi simplement séries de Fourier, et d'autre part les fonctions continues sommes de séries trigonométriques absolument convergentes, qui constituent l'algèbre de Wiener. Le point essentiel, mis en évidence par Wiener, est que ces espaces~$L^1$ et~$\ell^1$ sont des algèbres pour la convolution. Par transformation de Fourier, convolution et multiplication s'échangent. La convolution se trouve partout en analyse, mais en tant que notion c'est l'une des dernières nées: le terme de convolution ne s'est imposé que dans la seconde moitié du 20\ieme siècle. Au moins Wiener en avait fait sous le nom de Faltung un merveilleux usage, qui a servi de modèle à plusieurs égards à la théorie des distributions de Schwartz.

L'algèbre de Wiener a été un grand sujet de l'analyse au cours des années 1950 et 1960, et j'ai collaboré à cette époque à la découverte de ses secrets. Puis d'autres sujets sont passés au premier plan. Il est remarquable pour moi que les problèmes actuels de théorie du signal et de traitement des données donnent un nouveau départ à ce secteur de l'analyse de Fourier.\enlargethispage{\baselineskip}%

Wiener est connu à bien des titres. Le pense ici au mouvement brownien et aux séries de Fourier à coefficients aléatoires, un autre sujet sur lequel j'ai travaillé dans la tradition de Paley, Wiener et Zygmund. Le théorème de Candès, Romberg et Tao par lequel j'ai débuté cet exposé marie la transformation de Fourier discrète, l'espace~$\ell^1$ et les probabilités. C'est donc de plusieurs manières un héritage de Wiener.

Il est bon de conclure en revenant à ce théorème pour en indiquer le contexte et d'éventuels prolongements.

L'article \cite{CaRoTa} donne en introduction l'exposé d'une expérience numérique qui donne d'excellents résultats et qui viole l'usage établi par le théorème d'échantillonnage de Shannon. Il s'agit d'une figure appelée le \og fantôme de Logan-Shepp\fg, constituée par quelques plages en forme d'ellipses qui empiètent les unes sur les autres, et qui est un test classique pour les méthodes d'analyse du signal. On regarde sa transformée de Fourier sur 22 demi-droites issues de l'origine et régulièrement réparties, puis les valeurs qu'elle prend sur une vingtaine de points sur chaque demi-droite, au total 512 points.

On reconstitue exactement la figure par un procédé d'analyse convexe à partir des valeurs prises par la transformée de Fourier sur ces points, qui sont espacés les uns des autres d'au moins $\pi/22$, alors que le théorème d'échantillonnage exigerait pour être appliqué que l'espacement des points voisins soit $2 \pi / N$ avec $N =1000$ environ. Cet exemple n'est pas encore parfaitement compris, mais il repose sur la recherche du minimum dans l'algèbre de Wiener des normes des fonctions qui extrapolent la restriction à un ensemble convenable de la transformée de Fourier du signal considéré. Le cadre des groupes cycliques d'ordre $N$ pour~$G$ et $\hat G$ se prête bien à l'étude théorique, et l'échantillonnage aléatoire aboutit à l'énoncé que j'ai indiqué.

L'échantillonnage aléatoire est essentiel dans la démonstration, et il est également très intéressant en pratique. En expérimentant son utilisation on constate dans des cas simples que le choix de $\Omega$, l'effectif de l'échantillon, peut être rendu notablement plus petit que la valeur fixée par le théorème. L'article signale les résultats d'expériences allant dans ce sens.\enlargethispage{2\baselineskip}%

Quel est le champ d'application de ces méthodes ? Au contraire de la méthode des moindres carrés, qui consistait à gérer un ensemble de données surabondantes, il s'agit ici de restreindre les observations au minimum indispensable, un peu plus de $T$ si l'on sait que le signal est constitué de $T$ pics au plus.

L'utilisation de l'analyse convexe, qui rejoint ce que Fourier appelait l'analyse indéterminée, c'est à dire l'étude des inégalités et leur interprétation géométrique, s'est introduit progressivement dans ce domaine. Les normes~$\ell^1$ se sont avérées d'usage commode, lorsque la nature des problèmes suggérait la norme~$\ell^0$, c'est-à-dire le nombre de termes non nuls. Un exemple typique se trouve dans l'article de Donoho et Huo de 2001 dont la motivation est la reconstruction à partir d'un développement mixte en ondelettes et en Fourier, ou en ondelettes et en bandelettes, d'une image comprenant des plages et des contours. Le cadre général est celui déjà dressé en 1993 par Mallat et Zhang: on dispose d'un \og dictionnaire\fg constitué de plusieurs bases orthonormales, et il s'agit de l'utiliser pour obtenir de meilleurs résultats qu'avec une seule base. Dans l'article \cite{DoHu}, l'étude théorique de départ concerne les groupes cycliques d'ordre $N$,~$G$ et $\hat G$. Toute fonction définie sur~$G$ peut se traduire par Fourier en une fonction définie sur $\hat G$, et réciproquement. Il y a donc plusieurs représentations possibles pour la même fonction, en exprimant une partie sur~$G$ et une partie sur $\hat G$. Une représentation économique utilise le minimum de points de~$G$ et de $\hat G$, soit $T$ points de~$G$ et $\Omega$ points de $\hat G$ ; sa norme dans~$\ell^0$ est alors $T+\Omega$ ; et on peut espérer obtenir la meilleure représentation en minimisant la norme~$\ell^0$. C'est bien le cas lorsque $T+\Omega < \sqrt N$, ce ne l'est plus lorsqu'il y a égalité (description de l'exemple ? c'est la fonction indicatrice du sous-groupe d'ordre $n$ quand $N = n^2$, qui, à un facteur multiplicatif près, est sa propre transformée de Fourier: le minimum de la norme~$\ell^0$ est atteint de deux manières). Mais chercher le minimum dans~$\ell^0$ est impraticable. Le minimum dans~$\ell^1$ est praticable et donne le bon résultat quand $T+\Omega <\frac12\sqrt N$.

\part{Restrictions et prolongements dans~l'algèbre~de~Wiener}

\section*{Introduction}

Les travaux sur l'échantillonnage parcimonieux (compressed sensing) en théorie du signal ont mis en relief un très intéressant problème d'analyse de Fourier. En me référant aux notations d'Emmanuel Candès dans \cite{Can}, je le présenterai de la manière suivante.\enlargethispage{\baselineskip}%

$G$ est un groupe abélien localement compact et $\hat G$ son groupe dual. Soit $(S)$ une partie de $L^1(G)$, et $K$ une partie de $\hat G$. A toute $x\in L^1(G)$ on associe sa transformée de Fourier $\hat x$, qui appartient à l'algèbre de Wiener $A(\hat{G})$ ($=\cF L^1(G)$), et la restriction $\hat x$ à $K$, $ \hat{x}| _K$. Les $y\in L^1(G)$ tels que $ \hat{y}|_K = \hat{x}|_K$ constituent un convexe fermé dans $L^1(G)$, et les $\hat y$, prolongements de $ \hat{x}|_K$ dans $A( \hat{G})$, constituent un convexe fermé dans $A( \hat{G})$. Il s'agit de donner une condition sur $(S)$ et $K$ garantissant que, si $x\in (S)$, $\hat x$ est le prolongement minimal unique de $ \hat{x}|_K$ dans $A( \hat{G})$:
$$
\| \hat{x}_0\|_{A( \hat{G})} = \inf \{ \| \hat{y}\|_{A( \hat{G})} \mid \hat{y}|_K = \hat{x}|_K \} \Longleftrightarrow \hat{x}_0 = \hat{x}. \leqno(1)
$$

Dans \cite{Can} et dans les travaux qui précèdent, de Candès, Romberg et Tao \cite{CaRoTa} et antérieurement, sur un sujet voisin, de Donoho et Huo \cite{DoHu}, $G$ et $\hat G$ sont le groupe cyclique $\bbZ/N\bbZ$, $N$ étant un entier très grand. Le groupe $G$ est une représentation discrète des temps et le groupe $\hat G$ des fréquences. En désignant par $S$ le support du signal $x$, $(S)$ est traduit par une condition sur $S$, en l'occurrence dans \cite{Can} et \cite{CaRoTa} que le cardinal de $S$, $|S|$, est majoré par un nombre~$T$ petit par rapport~à~$N$:
$$
|S| = \|x\|_{\ell^o} \le T \ll N
\leqno(2)
$$
(la notation $\|\ \|_{\ell^o}$ a justement (2) pour définition). $K$ est un ensemble de fréquences, de cardinal $|K|$. Le théorème principal de \cite{Can} et \cite{CaRoTa} est que, si l'on fixe $|K|$ et que l'on choisit $K$ aléatoirement avec la probabilité naturelle, on retrouve le signal $x$ par le procédé (1) avec une probabilité voisine de 1 si $|K|$ est assez grand par rapport à $T \log N$, à savoir~que
$$
|K| \ge C\, T \log N
\leqno(3)
$$
avec $C\ge 23 (M+1)$ entraîne 
$$
P((1)) \ge 1-c N^{-M},
\leqno(4)
$$
où $c$ ne dépend que du signal.

Ainsi, si l'on sait que le signal $x$ est porté par $T$ points, avec $T\ll N$, on est loin d'avoir besoin des $N$ coefficients de Fourier de $x$ pour le reconstituer: en prenant $C=100$ dans (3), $|K|$ coefficients choisis au hasard suffisent pratiquement.

Les titre de \cite{CaRoTa} et \cite{DoHu} mettent en avant ce que les auteurs appellent \og uncertainty principles\fg, c'est-à-dire le fait (incontestablement certain) que les supports de $x$ et de $\hat x$ ne peuvent pas être petits simultanément. Les relations de Heisenberg, d'où est tiré le vocable de \og principe d'incertitude\fg, expriment simplement que l'on ne peut pas localiser simultanément un signal en temps et en fréquence. En physique, elles signifient que la représentation des particules par des points matériels est inadéquate à l'échelle quantique. En mathématiques, elles constituent un théorème de localisation sur transformées de Fourier. Pour $G=\bbZ/N\bbZ$, $\hat G$ est une autre copie de $\bbZ/N\bbZ$, et le théorème de localisation s'exprime ainsi:
$$
\|x\|_{\ell^o} \cdot\| \hat{x}\|_{\ell^o} \ge N.
\leqno(5)
$$
La démonstration, que j'indique en appendice, donne des informations supplémentaires en fonction des propriétés arithmétiques de $N$. On exprime souvent le théorème de localisation dans les $\bbZ/N\bbZ$ par une conséquence de (5), à savoir
$$
\|x\|_{\ell^o} +\| \hat{x}\|_{\ell^o} \ge 2 \sqrt{N}\,.
\leqno(6)
$$
L'égalité dans (6) a lieu lorsque $N$ est un carré parfait et que les supports de~$x$ et de $\hat x$ sont des translatés du sous-groupe cyclique d'ordre $\sqrt{N}$ du groupe $\bbZ/N\bbZ$, en particulier quand $x$ et $\hat x$ sont le \og peigne de Dirac\fg\ d'ordre $\sqrt{N}$, divisé par $\sqrt{N}$.

Dans \cite{DoHu}, Donoho et Huo tirent parti de cet exemple comme limitation naturelle pour la question qu'ils posent de la représentation d'un signal la plus économique possible en temps et en fréquences, c'est-à-dire l'expression de $x$, fonction définie sur $G=\bbZ/N\bbZ$, sous la forme
$$
x=y+z, \quad \|y\|_{\ell^o} + \| \hat{z}\|_{\ell^o}\ \text{minimum}.
\leqno(7)
$$
Quand $N$ est un carré parfait et que
$$
\|x\|_{\ell^o} = \| \hat{x}\|_{\ell^o} = \sqrt{N},
\leqno(8)
$$
on peut choisir $x=y$ ou $x=z$ dans (7): la représentation la plus économique n'est pas définie. Donoho et Huo montrent qu'elle est bien définie lorsque $x$ est de la forme
$$
x=y+z, \quad \|y\|_{\ell^o} + \| \hat{z}\|_{\ell^o} < \sqrt{N}
\leqno(9)
$$
et que, sous l'hypothèse additionnelle
$$
x=y+z, \quad \|y\|_{\ell^o} + \| \hat{z}\|_{\ell^o} < \frac{1}{2} \sqrt{N},
\leqno(10)
$$
on obtient la solution (7) en remplaçant la condition
$$
\|y\|_{\ell^o} + \| \hat{z}\|_{\ell^o} \ \text{minimum},
\leqno(11)
$$
pratiquement ingérable, par la condition
$$
\|y\|_{\ell^1} + \| \hat{z}\|_{\ell^1}\ \text{minimum},
\leqno(12)
$$
qui se prête au calcul.

Candès, Romberg et Tao utilisent leur théorème sur le prolongement minimal (formules (3) et (4)) pour proposer une amélioration considérable du théorème de localisation (6) sous la forme d'une proposition hautement probable: si l'on choisit au hasard une partie $S$ de $G=\bbZ/N\bbZ$ et une partie $S'$ de $\hat G$ telles~que
$$
|S| +|S'| < \gamma\, \frac{N}{\sqrt{\log N}}\,,
\leqno(13)
$$
alors, avec une probabilité aussi voisine de 1 que l'on veut par choix de $\gamma$, il n'existe aucune fonction $x$ définie sur $G$ à support dans $S$ telle que $\hat x$ ait son support dans~$S'$.

Par cet aperçu très partiel des travaux de Donoho et Huo de 2001 et de Candès, Romberg et Tao de 2006, on voit la portée en théorie du signal de la méthode de minimisation dans $\ell^1$. Cette méthode n'est pas nouvelle, et son histoire est retracée brièvement dans \cite{CaRoTa} p\ptbl494. Les initiateurs, au cours des années 1980, semblent avoir été des géophysiciens. La minimisation des normes~$\ell^1$ est une question d'analyse convexe dont on peut trouver les prolégomènes dans ce que Fourier appelait l'analyse indéterminée \cite{Fou}. La réduction de ce sujet d'analyse convexe à des techniques de programmation linéaire a été exposée en 1998 par Chen, Donoho et Saunders dans \cite{ChDoSa}. La minimisation de normes~$\ell^1$ est aujourd'hui un outil important de la théorie du signal.

Je vais m'attacher aux conditions de validité de (1) dans un contexte très différent: la reconstruction d'une fonction à spectre lacunaire à partir de sa restriction à un intervalle par prolongement minimal dans l'algèbre de Wiener.

\section*{Reconstruction des fonctions à spectre lacunaire}

Dans cette section on s'écarte de l'analyse de Fourier discrète. Le groupe $G$ sera $\bbZ^\nu$ $(\nu=1,2,\ldots)$ et $\hat G$ sera $\bbT^\nu$. On s'intéresse aux fonctions $ \hat{x} \in A(\bbT^\nu)$,
$$
\hat{x}(t) = \sum_{n\in S} x(n) e(n\cdot t), \quad \| \hat{x}\|_A = \sum_{n} |x(n)| < \infty
\leqno(14)
$$
($t=(t_1,\ldots, t_\nu)$, $n=(n_1,\ldots, n_\nu)$, $n\cdot t=n_1t_1 + \cdots n_\nu t_\nu$, $e(u)=\exp (2\pi iu)$) dont le spectre $S$ est lacunaire au sens que deux points distincts de $S$ sont toujours à une distance~$\ge d$:
$$
\inf_{\substack{s\in S\\s'\in S\\ s'\not= s}} |s'-s| \ge d
\leqno(15)
$$
$|\cdot|$ représente la distance euclidienne, et le premier membre de (15) s'appelle le \og pas\fg\ de $S$. Ici $K$ sera une boule fermée de rayon $r$ dans $\bbT^\nu$ $(0<r<\sfrac{1}{2})$. Le problème posé dans l'introduction est de donner une condition sur $d$ et $r$ ($\nu$ étant donné) garantissant que $\hat x$ est le prolongement minimal unique dans $A(\bbT^\nu)$ de sa restriction à $K$, $ \hat{x}|_K$, c'est-à-dire que la formule (1) est valable.

Pour faciliter la compréhension, limitons-nous d'abord au cas $\nu=1$. Alors~$K$ est un intervalle, et comme la solution du problème est invariante par les translations sur $\bbT$, on peut choisir $K=[-r,r]$. Ecrivons (14) sous la forme
$$
\hat{x} = \sum_{n\in S} x_n e_n, \quad \| \hat{x}\|_A = \sum_n |x_n| < \infty.
\leqno(16)
$$
(1) signifie que, quelle que soit $ \hat{z} \in A(\bbT)$ nulle sur $K$ et $\not\equiv 0$, soit
$$
\hat{z} = \sum_{n\in \bbZ} z_n e_n, \quad 0 < \sum_n |z_n| < \infty,\quad z|_K =0,
\leqno(17)
$$
on a
$$
\sum_{n\in S} |x_n+z_n| + \sum_{n\notin S} |z_n| > \sum_{n\in S} |x_n|.
\leqno(18)
$$
Dire que (18) a lieu pour toute suite $(x_n) \in \ell^1(S)$, c'est dire que
$$
\sum_{n \notin S} |z_n| > \sum_{n\in S} |z_n|
\leqno(19)
$$
(on a utilisé $\sup_{x_n}(|x_n| -|x_n+z_n|)=|z_n|$ ; c'est le même calcul qu'en \cite{CaRoTa} p\ptbl494). On cherche donc une condition sur $d$ et $r$ pour que (17) implique~(19).

Il y a bien des façons d'aborder le problème. Voici celle qui m'a donné les meilleurs résultats (loin sans doute d'être optimaux). On multiplie $\hat z$, qui est nulle sur $K$, par une fonction $\hat \varphi$ qui est presque entièrement portée par $K$. Le résultat est voisin de $O$, et il en est de même de sa transformée de Fourier, convolution de $z$ et de $\varphi$. Si $\varphi$ est presque entièrement portée par l'intervalle $[-d,d]$, la quasi-nullité de la convolution $z*\varphi$ est capable d'entraîner~(19).

Il est commode de partir de la fonction de Gauss
$$
\gamma(t) = e^{-\pi t^2} \qquad (t \in \bbR),
\leqno(20)
$$
de poser
$$
\gamma_a(t) = a \gamma(at) \qquad (a\gg\sfrac{1}{r})
\leqno(21)
$$
et de réduire $\gamma_a$ modulo 1 sous la forme
$$
\gamma_a^*(t) = \sum_{m\in \bbZ} \gamma_a(t-m) = \sum_{n\in \bbZ} \gamma (\sfrac{n}{a})\,e_n(t).
\leqno(22)
$$
C'est $\gamma_a^*$ qui nous servira de fonction $\hat \varphi$. Comme elle est paire, et décroissante et convexe sur l'intervalle $[r,\frac{1}{2}]$, sa restriction à $\bbT \ba K$ a une norme $\le 2\gamma_a^*(r)$ dans l'algèbre quotient $A(\bbT \ba K)$, d'où résulte que
$$
\| \hat{z} \gamma_a^*\|_{A(\bbT)} \le 2 \| \hat{z}\|_{A(\bbT)} \gamma_a^*(r).
\leqno(23)
$$

Pour chaque $s\in S$, intégrons sur $\bbT$ le produit $ \hat{z} \gamma_a^* e_{-s}$:
$$
\int \hat{z} \gamma_a^* e_{-s} = \sum_n z_n \gamma\big(\frac{n-s}{a}\big) = z_s + \sum_{n\not= s}
z_n \gamma \big(\frac{n-s}{a}\big)
\leqno(24)
$$
Multiplions par $\varepsilon_s = \sfrac{|z_s|}{z_s}$ et faisons la somme sur~$S$:
$$
\sum_s \varepsilon_s \int \hat{z} \gamma_a^* e_{-s} = \sum_{s\in S} |z_s| + \sum_{n\in \bbZ, s\in S,n\not= s} \varepsilon_s z_n \gamma(\frac{n-s}{a}).
\leqno(25)
$$
Dans cette dernière somme, distinguons les $n\in S$, que nous appellerons $s'$, et les $n\notin S$, et remarquons que le premier membre de (25) est majoré en valeur absolue par celui de (23). On obtient
$$
\Big| \sum_{s\in S} |z_s| + \sum_{s'\in S} A(s') z_{s'} + \sum_{n\notin S} B(n) z_n \Big| \le 2 \| \hat{z}\|_{A(\bbT)} \gamma_a^*(r)
\leqno(26)
$$
avec
$$
\begin{cases}
A(s') = \displaystyle\sum_{s\in S \ba \{s'\}} \varepsilon_s \gamma\big(\dfrac{s'-s}{a}\big)\\[15pt]
B(n) = \displaystyle\sum_{s\in S} \varepsilon_s \gamma\big(\dfrac{n-s}{a}\big)
\end{cases}
\leqno(27) 
$$
et (26) entraîne
$$
\sum_{n\notin S}
 |z_n| \big(|B(n)| + 2 \gamma_a^*(r)\big) \ge \sum_{s\in S} |z_s| \big(1- |A(s)| -2 \gamma_a^* (r)\big).
 \leqno(28)
$$
En vue d'obtenir (19), supposons
$$
\frac{1}{r} < a <d \quad \mathrm{et}\ a \ge 10
\leqno(29)
$$
et majorons $|B(n)|$, $\gamma_a^*(r)$ et $|A(s)|$.

Quand $n\notin S$, soit $h$ la distance de $n$ à $S$. Alors
\[
\tag*{(30)}
\begin{split}|B(n)| &\le \gamma \big(\sfrac{h}{a}\big) + \gamma \big(\sfrac{(d-h)}{a}\big) + 2 \sum_{j=1}^\infty \gamma \big(\sfrac{jd}{a}\big)\\
&\le \gamma \big(\sfrac{1}{a}\big) + \gamma \big(\sfrac{d}{2a}\big) + 2 \sum_{j=1}^\infty \gamma \big(\sfrac{jd}{a}\big)\\
&\le 1- \big(\sfrac{3\pi}{4a^2}\big) + 2 \exp \big(-\pi\sfrac{d^2}{4a^2}\big).
\end{split}
\]
Pour $s\in S$,
$$
|A(s)| \le 2 \sum_{j=1}^\infty \gamma \big(\sfrac{jd}{a}\big) \le 2\exp \big( -\pi \sfrac{d^2}{4a^2}\big).
\leqno(31)
$$
Enfin
$$
2\gamma_a^*(r) \le 3a \exp (-\pi a^2 r^2).
\leqno(32)
$$
Supposons
$$
r > \frac{2}{a} \, \sqrt{\log a}\,; 
\leqno(33)
$$
alors
$$
2\gamma_a^*(r) \le \frac{\pi}{8a^2}.
\leqno(34)
$$
Supposons
$$
d > 2a \sqrt{\log a}\,; 
\leqno(35)
$$
alors
$$
\begin{cases}
|A(s)| \le \dfrac{\pi}{4a^2}\\
|B(n)| \le 1- \dfrac{2\pi}{4a^2}
\end{cases}
\leqno(36)
$$
et (19) est bien établie.

L'une des conditions
$$
d> \frac{5}{r} \log \frac{1}{r} \qquad \big(r<\frac{1}{10}\big)
\leqno(37)
$$
ou
$$
r > \frac{5}{d} \log d \qquad (d>10)
\leqno(38)
$$
entraîne donc, en choisissant $a=\sqrt{\sfrac{d}{r}}$, ce qui implique (29), que pour toutes les fonctions $\hat z$ vérifiant (17) on~a~(19).

L'étude faite en dimension 1 se traduit sans peine en dimension $\nu>1$, au prix d'une augmentation des constantes.

On prend pour $K$ la boule $|t|<r$ et on remplace $\gamma$, $\gamma_a$ et $\gamma_a^*$~par
$$
\gamma_\nu(t) = e^{-\pi|t|^2} \qquad (t\in \bbR^\nu)
\leqno(39)
$$
$$
\gamma_{a,\nu}(t) = a^\nu\gamma(at) \qquad \big(a\gg \frac{1}{r}\big)
\leqno(40)
$$
\[
\begin{split}
\gamma_{a,\nu}^*(t) &=\sum_{m\in \bbZ^\nu} \gamma_{a,\nu} (t-m)\\
&=\sum_{n\in \bbZ^\nu} \gamma_\nu \big(\sfrac{n}{a}\big)\,e(n\cdot t)
\end{split}
\tag*{(41)}
\]

La formule (28) reste valable en remplaçant $2\gamma_a^*(r)$ par $\|\gamma_{a,\nu}^*\|_{A(\bbT^\nu\ba K)}$, et~$A(s)$ et $B(n)$ par~$A_\nu(s)$ et $B_\nu(n)$, définis comme en (27) avec $\gamma_\nu$ au lieu de~$\gamma$. La majoration de $\|\gamma_{a,\nu}^*\|_{A(\bbT^\nu\ba K)}$ peut se faire ici en estimant les dérivées d'ordre $\le \nu$ de $\gamma_{a,\nu}^*$ sur $\bbT^\nu\ba K$, et au lieu de (32) on a
$$
\|\gamma_{a,\nu}^*\|_{A(\bbT^\nu\ba K)} \le C_\nu a^\nu \exp (-\pi a^2 r^2),\quad C_\nu <\nu!
\leqno(42)
$$
La majoration de $|B_\nu(n)|$ se fait comme en (30), en distinguant le point $s$ le plus voisin de $n$ (ou l'un de ces points s'il y en a plusieurs) et les points de $S$ dans la boule de centre $s$ et de rayon $\frac{d}{2}$, dont le nombres est $\le 2^\nu$. Celle de $|A_\nu(s)|$ est analogue, en plus simple. On obtient
$$
\left\{
\begin{array}{ll}
|B(n)| &\le 1-\dfrac{3\pi}{4a^2} + 2^{\nu+1}
 \exp \big(-\pi \dfrac{d^2}{4a^2}\big)\\
 |A(s)| &\le 2^{\nu+1} \exp \big(-\pi \dfrac{d^2}{4a^2}\big),
\end{array}
\right.
\leqno(43)
$$
d'où finalement l'énoncé suivant.

\begin{theoreme*}
Soit $ \hat{G}=\bbT^\nu(\nu=1,2,\ldots)$. Pour que toute fonction $ \hat{x}\in A( \hat{G})$, dont le pas du spectre est $d$ ($d\ge 10$) s'obtienne à partir de sa restriction $ \hat{x}|_K$ à une boule $K$ de rayon $r$ comme prolongement minimal de $ \hat{x}|_K$ dans $A( \hat{G})$, il suffit que l'une des conditions suivantes soit réalisée:
\begin{gather*}
\tag{44}
r > \frac{5}{d} \ \sqrt{\nu}\, \log(\nu d)\\
\tag{45}
d > \frac{5}{r} \ \sqrt{\log(\sfrac{1}{r}) +\nu-1}
\end{gather*}
\end{theoreme*}

(Pour $\nu=1$, ce sont les formules (38 et (37)).

La constante 5 peut être abaissée quand $d$ est grand, et doit être ajustée pour couvrir tous les cas.

\section*{Discussion et commentaires}

Je commencerai par le cas $\nu=1$.

\vskip2mm

\textbf{1.} On peut transcrire le théorème dans le contexte donné dans l'introduction, où $G$ et $\hat G$ sont le groupe cyclique $\bbZ/N\bbZ$, $N$ très grand, $G$ étant une représentation discrète des temps et $\hat G$ des fréquences. Alors $x$ est un signal porté par des temps séparés par des intervalles de longueurs $\ge d$, et $K$ est une bande de fréquences de largeur $\ge 2rN$ (si l'on représente $\hat G$ par les racines~$N$\nobreakdash-ièmes de l'unité, $K$ est la partie de $\hat G$ contenue dans un arc du cercle unité). On obtient l'énoncé que voici:

Si le support $S$ du signal $x$ (porté par $\bbZ/N\bbZ$) a un pas $\ge d$, c'est-à-dire~si
$$
\inf_{\substack{s\in S\\s'\in S\\ s'\not=s}} |s'-s| \ge d \qquad (d\ge 10)
\leqno(46)
$$
on peut obtenir $x$ à partir de $|K|$ fréquences consécutives (dans $\bbZ/N\bbZ$) par le procédé indiqué en (1), dès que
$$
|K| \ge 10\, \frac{N}{d} \, \sqrt{\log d}
\leqno(47)
$$

Le défaut de cet énoncé par rapport à celui de Candès, Romberg et Tao est dans l'hypothèse faite sur $S$: (46) implique $|S| \le \sfrac{N}{d}$, mais $|S| \le T$ n'implique rien sur $d$. Cependant, si $T$ est donné et si $S$ s'obtient en disposant~$T$ points au hasard sur $\bbZ/N\bbZ$, on a $\sfrac{N}{d} > \sfrac{10}{T^2}$ avec une probabilité $>0,95$, donc il est probable que
$$
|K| \ge 100\, T^2 \sqrt{\log T}
\leqno(48)
$$
suffit pour que (1) soit valable. L'avantage est qu'il s'agit de fréquences consécutives, c'est-à-dire d'une bande de fréquences.

Cela serait à comparer aux résultats connus.

\vskip2mm

\textbf{2.} On a beaucoup étudié les séries trigonométriques lacunaires et les fonctions qu'elles représentent:
$$
\hat{x}(t) = \sum_{n\in S} x_n e_n(t) = \sum_{n\in S} x_n e^{2\pi int} \qquad (t\in \bbT).
\leqno(49)
$$C'est de ces séries et fonctions que nous nous sommes occupés à partir de (16), sous la condition $\sum |x_n| < \infty$. On peut les étudier sous d'autres conditions, en particulier $\sum |x_n|^2 < \infty$.

En gros, si $S$ est lacunaire, beaucoup de propriétés de $\hat x$ se lisent sur les restrictions $ \hat{x}|_K$, où $K$ est un ensemble convenable. Je me limiterai ici au cas des intervalles.

\emph{Première propriété}: la nullité. Quand est-il vrai, sous l'hypothèse (49), que $ \hat{x}|_K=0$ entraîne $ \hat{x}=0$ ? Autrement dit, que $\hat x$ soit bien défini par les valeurs qu'elle prend sur $K$ ? Cela se traduit lorsque $K=[-r,r]$ par une question sur les fonctions entières de type exponentiel $\le \pi(1-r)$ bornées (ou $L^1$, ou~$L^2$...) sur la droite réelle, $X(z)$ $(z\in \bbC)$: à quelle condition leur nullité sur $\bbZ \ba S$ entraîne-t-elle leur nullité ? Plus simplement, $S$ étant donnée, quelle est la borne inférieure des $r$ tels qu'il en soit ainsi ? La réponse est fournie par la théorie de Beurling et Malliavin \cite{BeMa,Koo}. Associons à toute suite $\Lambda \subset \bbZ$ sa fonction de décompte $n(t)$, dont la variation sur un intervalle égale de nombre de points de $\Lambda$ sur cet intervalle, et convenons de dire que $\Lambda$ admet $\Delta$ comme densité-BM si
$$
\int_\bbR \frac{|n(t)-\Delta t|}{1+t^2} dt < \infty.
\leqno(50)
$$
La densité-BM intérieure de $S$ se définit naturellement comme la borne supérieure des densités-BM des suites $\Lambda$ contenues dans $S$. C'est justement la borne inférieure dans $|K|$ (longueur de $K$) tels que $ \hat{x}|_K=0$ entraîne $ \hat{x}=0$ $(K=[-r,r])$.

En particulier, pour que $\hat x$ soit bien définie par les valeurs qu'elle prend sur n'importe quel intervalle, si petit soit-il, il faut et il suffit que $S$ ne contienne aucune suite $\Lambda$ ayant une densité-BM $>0$. Il n'est pas nécessaire que $S$ soit très dispersée ; il suffit qu'elle ait de grandes lacunes. 

\vskip2mm

\emph{Seconde propriété}: l'appartenance à l'algèbre de Wiener. Rappelons que $A(\bbT) = \cF \ell^1(\bbZ)$ et que $A(K)$ est l'algèbre quotient de $A(\bbT)$ par l'idéal constitué par les fonctions nulles sur $K$, autrement dit, l'algèbre des restrictions à $K$ des fonctions $\in A(\bbT)$. On considère les fonctions $\hat x$ de (49) avec $x\in \ell^1(\bbZ)$. Quelle est la borne inférieure des longueurs des intervalles $K$ tels que $\hat {x}|_K \in A(K)$ implique $ \hat{x} \in A (\bbT)$ ? Réponse: c'est la densité uniforme extérieure de $S$, c'est-à-dire la borne inférieure des densités $D$ des suites $\Lambda$ telles que
$$
n(t) =Dt+O(1) \qquad (t \to \pm \infty) \leqno(51)
$$
et qui contiennent $S$ ; on exprime (51) en disant que $\Lambda$ a $D$ pour densité uniforme~\cite{Kah2}.

\emph{Variante}: l'appartenance aux espaces $A^{\alpha,r}
(\bbT) = \cF \ell^{\alpha,r}
(\bbZ)$. La définition de $\ell^{\alpha,r}
(\bbZ)$ ($1\le \alpha\le \infty,\ r\in \bbR$) est
\[\tag*{(52)}
\ell^{\alpha,r}
(\bbZ) =
\begin{cases}
\Big\{ x \mid \sum_{n\in \bbZ} |x_n|^\alpha (1+|n|)^{r\alpha} < \infty \Big\}& \text{si } 1\le\alpha< \infty,\\[5pt]
\Big\{ x \mid\sup_{n\in \bbZ}
 |x_n| (1+|n|)^r < \infty \Big\}& \text{si } \alpha=\infty.
\end{cases}
\]
En général, $A^{\alpha,r}(\bbT)$ est un espace de distributions de Schwartz, dont on sait définir les restrictions à des intervalles ouverts $J$ ; ces restrictions constituent l'espace $A^{\alpha,r}(J)$. Question: quelle est la borne inférieure des $r$ tels que $ \hat{x}|_J \in A^{\alpha,r}(J)$ $(J=]-r,r[)$ entraîne $ \hat{x} \in A^{\alpha,r}(\bbT)$ ? Réponse: elle est la même que précédemment, indépendante de $\alpha$~et~$r$.

Pour $r=0$ et $\alpha=1$, on retrouve $A(\bbT)$. Pour $r=0$ et $\alpha=2$, on trouve $L^2(\bbT)$, et c'est l'origine de la théorie. Sans se borner à des $\Lambda\subset \bbZ$ mais en considérant des $\Lambda$ réelles \og régulières\fg\ au sens que leur pas
\[
\inf_{\substack{\lambda\in\Lambda\\\lambda'\in\Lambda\\\lambda\not=\lambda'}}|\lambda'-\lambda|
\]
 est $>0$, Paley et Wiener ont appelé pseudo-période attachée à $\Lambda$ la borne inférieure des $\ell>0$ tels que les normes dans $L^2(J)$ des fonctions
$$
f(t) = \sum x_\lambda\ e^{i\lambda t}
\leqno(53)
$$
soient équivalentes pour tous les intervalles $J$ de longueur $\ell$. Au facteur $2\pi$ près, la valeur de la pseudo-période est la densité uniforme extérieure de~$\Lambda$~\cite{Kah1}.

\vskip2mm

\emph{Troisième propriété}: la continuité. C'est la question la plus naturelle, mais la plus difficile. Je renvoie à l'étude d'Yves Meyer \cite{Meyer} pour les meilleurs résultats connus.

\vskip2mm

\emph{Quatrième propriété}: l'analycité. Sous une forme un peu différente, c'est la première question qui a été examinée (Fabry en 1896, P\'olya en 1942), mais la dernière ayant eu sa solution complètement explicitée \cite[II p\ptbl52]{Koo} \cite[p\ptbl199]{KaLe}. Quelle est la borne inférieure des longueurs $|J|$ telles que, si $\hat x$ est analytique sur~$J$, $\hat x$ est analytique partout sur $\bbT$ ? C'est la densité maximale de $S$, à~savoir la densité extérieure quand, au lieu de (50) et de (51), la densité $D$ d'une suite~$\Lambda$ est définie~par
$$
n(t) = Dt+o(t) \qquad (t\to \pm \infty).
\leqno(54)
$$
P\'olya a montré que la densité maximale est atteinte: il existe une suite contenant $S$ dont la densité (au sens de (54)) est la densité maximale de~$S$.

Lorsque le pas de $S$ est $\ge d$, toutes les densités extérieures de $S$ sont $\ge\sfrac{1}{d}$. On pourrait donc espérer reconstruire $\hat x$ à partir de $ \hat{x}|_K$ $(K=[-r,r])$ quand $|K| = 2r > \sfrac{1}{d}$. Nous avons vu (formule (44)) que c'est possible quand $|K| > \frac{10}{d} \sqrt{\log d}$. On peut améliorer le résultat en diminuant le second membre, mais pas au dessous de~$\sfrac{1}{d}$.

Aucune densité n'est adaptée à la reconstruction d'une fonction $\hat x$ appartenant à $A(\bbT)$ à partir d'une restriction $ \hat{x}|_K$ par prolongement minimal dans $A(\bbT)$. En effet, soit $ \hat{z} \in A(\bbT)$ une fonction nulle sur $K$. Soit $\hat x$ une somme partielle de sa série de Fourier telle que $\| \hat{z}- \hat{x}\|_{A(\bbT)} <\| \hat{x}\|_{A(\bbT)}$ ; alors le procédé (1) échoue pour $\hat x$, et le support $S$ de $x$ est fini, donc de densité nulle. L'hypothèse naturelle pour cette reconstruction est bien que la pas de $S$ soit convenablement minoré. Le théorème de Candès, Romberg et Tao a le grand mérite de prendre pour hypothèse une majoration de $|S|$ sans rien supposer de son pas ; mais il faut alors renoncer à prendre pour $K$ un intervalle, et choisir pour $K$ un ensemble fini aléatoire.

\vskip2mm

Quittons le cas $\nu=1$.

\vskip2mm

\textbf{3.} Pour $\nu \ge 2$, on est loin d'avoir des réponses aussi précises. Par exemple, on ne sait pas pour quelles valeurs de $r$ l'hypothèse
$$
\inf_{\substack{s\in S\\s'\in S\\ s'\not= s}} |s'-s| \ge d,
\leqno(55)
$$
$|\ |$ représentant la norme euclidienne, entraîne l'équivalence des normes\break $\| \hat{x}\|_{L^2(B)}$ pour toutes les boules $B$ de rayon $r$. La relation entre la géométrie de $S$ et celle des $K$ associés est un problème intéressant, amorcé dans \cite{Kah1}.

\section*{Appendice 1. Sur le principe de localisation de Tao}

$G$ et $\hat G$ étant deux groupes cycliques d'ordre $N$ en dualité, $x$ une fonction définie sur $G$ à valeurs complexes et $\hat x$ sa transformée de Fourier, il s'agit d'établir que
$$
\|x\|_{\ell^o} \cdot \| \hat{x}\|_{\ell^o} \ge N.
\leqno(56)
$$
Voici la preuve de Tao, ingénieuse et facile.

Soit $S$ le support de $x$, et $|S|=\|x\|_{\ell^o}=s \not= 0$. Rappelons que
$$
\hat{x}(t) = \sum_{n\in S} x_n e^{2\pi int/N}.
$$
Supposons $ \hat{x} (t_1) = \hat{x} (t_2) = \cdots = \hat{x}(t_s) =0$ avec $t_1,t_2,\ldots, t_s$ en progression arithmétique, tous distincts modulo $N$. Les $x_n$ $(n\in S)$ sont solutions d'un système linéaire homogène dont le déterminant, du type Van der Monde, est~\hbox{$\not= 0$}, donc ils sont tous nuls, contrairement à l'hypothèse. Ainsi, sous l'hypothèse $\|x\|_{\ell^o} = s>0$, $\hat x$ ne peut pas s'annuler en $s$ points formant une progression arithmétique.

Les intervalles d'entiers contigus au support $\hat S$ de $\hat x$ ont donc au plus $s$ éléments, ce qui entraîne $s| \hat{S}| \ge N$, c'est-à-dire (56).

On a égalité dans (56) lorsque $S$ et $\hat S$ sont des translatés de sous-groupes cycliques d'ordre $s$ et $\hat s$, de $G$ et $\hat G$, avec $s \hat{s}=N$.

\section*{Appendice 2. Reconstruction d'un signal par la méthode de~Candès, Romberg et Tao (variantes~déterministe~et~probabiliste)}
\let\oldtheequation\theequation
\renewcommand{\theequation}{\oldtheequation$'$}

Le signal est une fonction $x(t)$, $t\in G=\bbZ_N=\bbZ/N\bbZ$, portée par une partie~$S$ de $G$, de cardinal $T$. Sa transformée de Fourier est
\[
\wh x(t)=\frac{1}{\sqrt N}\sum_{t\in S}x(t)e(-\omega t/N),\quad\omega\in\wh G=\bbZ_N,\ e(u)=e^{2\pi i u}.
\]
On veut reconstituer le signal en n'utilisant qu'une partie $\Omega$ des fréquences $\omega$, par le procédé que voici:
\begin{enumeratei}
\item\label{enum:1}
$\wh x$ est le prolongement minimal de $\wh x_{|\Omega}$ dans $A(\wh G)$ ($=\cF L^1(G)$).

\indent Quand est-ce possible? On va se donner des conditions suffisantes, emboîtées. J'écrirai indifféremment $x(t)$ ou $x_t$. \eqref{enum:1} signifie que pour toute $\wh z$ nulle sur~$\Omega$ et $\not\equiv0$,
\[
\sum_{t\in S}|x_t+z_t|+\sum_{t\in G\ba S}|z_t|>\sum_{t\in S}|x_t|
\]
et il suffit pour cela que 
\[
\sum_{t\in G\ba S}|z_t|>\sum_{t\in S}|z_t|.
\]
Donc il suffit pour avoir \eqref{enum:1} que

\item\label{enum:2}
Pour toute $\wh z$ nulle sur $\Omega$ et $\not\equiv0$, $\sum_{t\in G\ba S}|z_t|>\sum_{t\in S}|z_t|$.

\indent
Pour traduire \eqref{enum:2}, introduisons $\lambda_t$ ($t\in G$) telle que $|\lambda_t|=1$ et $\lambda_t\ov z_t=|z_t|$. Soit $P(t)$ ($t\in G$) tel que $\wh P$ soit porté par $\Omega$. Alors
\[
\int_G P\ov z=\int_{\wh G}\wh P\ov{\wh z}=0,
\]
soit
\begin{equation}\label{eq:1}
\int_S+\int_{G\ba S}P\ov z=0.
\end{equation}
Supposons que 
\begin{equation}\label{eq:2}
\begin{cases}
|P(t)-\lambda_t|<1/2&\text{quand }t\in S,\\
|P(t)|<1/2&\text{quand }t\in G\ba S.
\end{cases}
\end{equation}
Alors
\[
\begin{cases}
\dpl\Big|\int_SP\ov z-\lambda\ov z\Big|<\frac12\int_S|\ov z|,&\quad\text{d'où }\dpl\Big|\int_SP\ov z\Big|>\frac12\int_S|\ov z|,\\[10pt]
\dpl\Big|\int_{G\ba S}P\ov z\Big|<\frac12\int_{G\ba S}|\ov z|
\end{cases}
\]
et, compte-tenu de \eqref{eq:1}, on a
\[
\int_{G\ba S}|\ov z|>\int_S|\ov z|,
\]
c'est-à-dire la conclusion de \eqref{enum:2}. Il suffit donc, pour avoir \eqref{enum:2}, que
\item\label{enum:3}
quels que soient les $\lambda_t$ ($t\in G$) de module $1$, il existe $P(t)$ ($t\in G$) tel que $\wh P$ soit porté par $\Omega$ et qu'on ait \eqref{eq:2}.

\indent
Reste à construire $P$. Pour cela, introduisons le noyau
\begin{equation}\label{eq:3}
K(t)=\sum_{\omega\in\Omega}e(\omega t/N)\quad(t\in G)
\end{equation}
et posons
\begin{align*}
P(t)&=\sum_{t'\in S}\lambda_{t'}\,\frac{K(t-t')}{K(0)}\\
&=\sum_{\omega\in\Omega}\frac{1}{K(0)}\sum_{t'\in S}\lambda_{t'}\,e(-\omega t'/N)\,e(\omega t/N).
\end{align*}
Ainsi, $\wh P$ est porté par $\Omega$. Supposons maintenant que

\item\label{enum:4}
le noyau $K(t)$ défini par \eqref{eq:3} vérifie
\[
|K(t)|<\frac{1}{2T}\,K(0)\quad\text{quand $t\neq0$ (rappelons que $T=|S|$).}
\]
Alors \eqref{eq:2} est vérifiée donc \eqref{enum:3} a lieu, donc \eqref{enum:2}, donc \eqref{enum:1}. La condition \eqref{enum:4} entraîne \eqref{enum:1} pour tous les signaux $x(t)$ portés par $T$ points.


Ce résultat facile est la variante déterministe du théorème CRT. Je comparerai plus tard la portée du théorème et des variantes.

La recherche des noyaux \eqref{eq:3} ayant la propriété \eqref{enum:4} est intéressante en elle-même. Les noyaux s'appellent idempotents parce qu'ils coïncident avec leur carré de convolution; la propriété \eqref{enum:4} traduit une certaine parenté entre~$\Omega$ et le groupe $\wh G$ entier. Pour aborder cette recherche, nous allons emprunter l'approche de Candès, Romberg et Tao.

Selon \cite{CaRoTa}, choisissons maintenant pour $\Omega$ la partie de $\wh G$ définie par une sélection aléatoire de paramètre $\tau$: pour chaque $n\in\wh G$, l'événement $\{n\in\Omega\}$ a pour probabilité $\tau$, et ces événements sont indépendants. En posant \[
X_n=\mathbf{1}_{n\in\Omega},
\]
les $X_n$ sont des variables de Bernoulli de paramètre $\tau$, indépendantes:
\[
\text{distribution de }X_n=B(1,\tau).
\]
Le cardinal $|\Omega|$ de $\Omega$ est $\sum X_n$, qui est voisin de $N\tau$:
\[
\text{distribution de }|\Omega|=B(N,\tau),
\]
et quand $|\Omega|$ est fixé, la distribution de l'ensemble aléatoire $\Omega$ est uniforme sur l'ensemble des parties de $\wh G$ à $|\Omega|$ éléments. Dans la suite, ce qui sera dit de l'ensemble aléatoire $\Omega$ sera également valable pour un ensemble de $[N\tau]$ éléments distribués au hasard sur $\wh G$ selon la probabilité naturelle.

Précisons ce point. Pratiquement, $N$ et $N\tau$ sont grands et on peut se représenter $B(N,\tau)$ comme une distribution normale $\cN(N\tau,N\tau(1-\tau))$ (espérance $N\tau$, variance $N\tau(1-\tau)$). Si, comme dans la suite, $\tau N=C(T)\log N$, on aura
\begin{equation}\label{eq:4}
P\big(\big||\Omega|-\tau N\big|>\lambda\sqrt{\log N}\big)=O(N^{-\frac12\lambda^2/C(T)})\quad(N\to\infty).
\end{equation}
Comme exercice, cherchons à démontrer cette formule. Soit $u>0$; on~a
\begin{equation}\label{eq:5}
E(e^{uX_n})=1-\tau+\tau e^u<\exp(\tau(e^u-1)).
\end{equation}
Si de plus $u<1$,
\begin{align*}
E(e^{uX_n})&<e^{\tau(u+u^2)},\\
E(e^{u\sum_1^NX_n})&<e^{N\tau(u+u^2)},\\
P\big(\textstyle\sum_1^NX_n-\tau N-\lambda\sqrt{\log N}>0\big)&<E\Big(\exp\big(u(\textstyle\sum_1^NX_n-\tau N-\lambda\sqrt{\log N})\big)\Big)\\
&<\exp(N\tau u^2-\lambda u\log N).
\end{align*}
Le minimum du second membre est atteint quand $u=(\lambda\log N)/2N\tau$ et vaut
\[
\exp\Big(-\frac{\lambda^2\log^2N}{4N\tau}\Big).
\]
En remplaçant $\sum_1^NX_n-\tau N$ par $-\sum_1^NX_n+\tau N$ et en ajoutant les probabilités obtenues, on a
\[
P\big(\big|\textstyle\sum_1^NX_n-\tau N\big|>\lambda\sqrt{\log N}\big)<2\exp\Big(-\dfrac{\lambda^2\log^2N}{4N\tau}\Big),
\]
c'est-à-dire, quand $\tau N=C(T)\log N$,
\begin{equation}\label{eq:6}
P\big(\big||\Omega|-\tau N\big|>\lambda\sqrt{\log N}\big)=2N^{-(\lambda^2\log^2N)/4N\tau}.
\end{equation}
L'exposant $1/4$ dans \eqref{eq:6} au lieu de $1/2$ dans \eqref{eq:4} tient à la majoration trop généreuse de $e^u-1$ par $u+u^2$.

Nous identifierons désormais $\wh G$ à $\bbZ_N$. Au lieu de \eqref{eq:3}, nous avons
\[
K(t)=\sum_{n\in\bbZ_N}X_ne(nt/N)
\]
et en particulier $K(0)=\sum X_n$. Étant donné $N$ et $T$, on va chercher à déterminer $\tau$ de façon que \eqref{enum:4} soit très probable.

Le contraire de \eqref{enum:4} est
\[
\exists t\in G\ba\{0\},\quad|K(t)|\geq\frac{1}{2T}\,K(0).
\]
Choisissons un entier $\nu\geq3$. Pour chaque $z\in\bbC\ba\{0\}$ il existe un $\varphi\in\{2k\pi/\nu\mid k=1,2,\dots,\nu\}$ tel que $\reel ze^{-i\varphi}\geq|z|\cos\pi/\nu$. Posons $a=\cos\pi/\nu$. Le contraire de \eqref{enum:4} implique
\[
\exists t,\varphi,\quad\reel K(t)e^{-i\varphi}\geq\frac{a}{2T}\,K(0).
\]
Fixons $t$ et $\varphi$, posons
\begin{align*}
Y&=\reel K(t)e^{-i\varphi}-\frac{a}{2T}\,K(0)\\
&=\sum_{n\in\bbZ_N}X_n\Big(\cos(2\pi nt/N-\varphi)-\frac{a}{2T}\Big)
\end{align*}
et cherchons à majorer $P(Y\geq0)$. On a
\[
P(Y\geq0)<E\exp uY\quad\text{pour tout $u>0$},
\]
et, utilisant \eqref{eq:5},
\begin{align*}
E\exp uY&=\prod_{n\in\bbZ_N}\exp\Big(u\big(\cos(2\pi nt/N-\varphi)-\frac{a}{2T}\big)X_n\Big)\\
&=\prod_{n\in\bbZ_N}\exp\tau\Big(\exp\Big(u\big(\cos(2\pi nt/N-\varphi)-\frac{a}{2T}\big)\Big)-1\Big)\\
&=\exp\tau\sum_{n\in\bbZ_N}\sum_{k=1}^\infty\frac{u^k}{k!}\Big(\cos(2\pi nt/N-\varphi)-\frac{a}{2T}\Big)^k\\
&=\exp\tau\sum_{k=1}^\infty\frac{u^k}{k!}\sum_{n\in\bbZ_N}\Big(\cos(2\pi nt/N-\varphi)-\frac{a}{2T}\Big)^k.
\end{align*}
Supposons  $N$  impair. Comme $t\neq0$, donc aussi $2t\neq0$, la dernière somme vaut
\[
\begin{cases}
-\dfrac{aN}{2T}&\text{quand }k=1,\\[10pt]
\dfrac N2+N\dfrac{a^2}{4T^2}&\text{quand }k=2
\end{cases}
\]
et elle est toujours majorée par $N2^k$, donc
\[
E\exp uY\leq\exp\tau N\Big[-\dfrac{a}{2T}\,u+\Big(\frac12+\frac{a^2}{4T^2}\Big)\frac{u^2}{2}+\sum_{k=3}^\infty\frac{(2u)^k}{k!}\Big].
\]
Pour minimiser la somme des deux premiers termes du crochet, prenons
\[
u=\frac{a}{2T}\Big(\frac12+\frac{a^2}{4T^2}\Big)^{-1}.
\]
On obtient
\begin{align*}
E\exp uY&\leq\exp\tau N\Big[-\dfrac{a^2}{8T^2}\Big(\frac12+\frac{a^2}{4T^2}\Big)^{-1}+\sum_{k=3}^\infty\frac{1}{k!}\Big(\frac{2a}{T}\Big)^k\Big]\\
&\leq\exp\Big[-\tau N\Big(\frac{a^2}{4T^2+2}-\frac{2}{T^3}\Big)\Big],
\end{align*}
donc
\[
P(\text{non}\eqref{enum:4})<N\nu\exp\Big[-\tau N\Big(\frac{a^2}{4T^2+2}-\frac{2}{T^3}\Big)\Big].
\]
Si l'on choisit
\[
\tau N=4C(T^2+1)\log N,\quad C=1+\delta>1,
\]
puis $\nu$ assez grand, donc $a$ assez voisin de $1$, et si $T$ est assez grand, cette probabilité est $O(N^{-\delta'})$ quand $N\to\infty$ pour tout $\delta'<\delta$.
\end{enumeratei}

En conclusion, si l'on choisit
\begin{equation}\label{eq:7}
|\Omega|\geq4(1+\delta)(T^2+1)\log N,\quad\delta>0,
\end{equation}
si $T$ est assez grand, la probabilité de \eqref{enum:4} est $1-O(N^{-\delta'})$ ($N\to\infty$) pour tout $\delta'<\delta$. Il en est de même pour la probabilité d'avoir \eqref{enum:1} pour \emph{toutes} les fonctions $x(t)$ portées par~$T$ points de $G$.

C'est la variante aléatoire que j'avais annoncée.

\medskip
Le théorème de Candès, Romberg et Tao dit que, étant donné \emph{une} fonction $x(t)$ portée par $T$ points de $G$, si l'on choisit
\begin{equation}\label{eq:8}
|\Omega|\geq22(1+\delta)T\log N,\quad\delta>0,
\end{equation}
la probabilité d'avoir \eqref{enum:1} pour \emph{cette} fonction est $1-O(N^{-\delta})$ (voir \cite[Th\ptbl2.1]{Can}; d'après \cite[(3.34)]{CaRoTa} on s'attendrait à $23$ au lieu de $22$ dans \eqref{eq:8}). Dans leurs énoncés, les auteurs n'insistent pas sur le fait que l'événement dont ils estiment la probabilité dépend du signal $x(t)$, mais il en est bien ainsi dans la démonstration.

L'intérêt de la version déterministe, outre la facilité de sa démonstration, est de montrer la portée possible en théorie du signal du problème de la concentration des idempotents dans $\bbZ_N$.

La condition \eqref{eq:7} de la version aléatoire est plus exigeante que \eqref{eq:8} dès que $T\geq6$. L'intérêt de la version aléatoire est qu'une fois $\Omega$ choisi, \eqref{enum:1} est valable pour tous les signaux portés par $T$ points; l'événement dont on estime la probabilité ne dépend que de $T$. Le prix à payer est de faire intervenir $T^2$ au lieu de $T$.

Pourrait-on remplacer $T^2$ par une quantité $o(T^2)$ ($T\to\infty$)? Si l'on s'attache à la probabilité de \eqref{enum:4}, la réponse est négative. En effet, \eqref{enum:4} entraîne
\[
N|\Omega|=\sum_{t\in G}|K(t)|^2\leq K(0)^2\Big(1+\frac{N-1}{4T^2}\Big)=|\Omega|^2\frac{N+4T^2-1}{4T^2},
\]
donc
\[
|\Omega|\geq4T^2\,\frac{N}{N+4T^2-1}.
\]

Peut-on retrouver et étendre le théorème CRT en adaptant la démonstration de la version aléatoire? La réponse est positive, mais il serait trop long de la donner ici.

\Subsection*{Exercices}

\begin{enumeratea}
\item
Montrer que, pour la valeur choisie pour $a$,
\[
\sum_{k=3}^\infty\frac{u^k}{k!}\sum_{n\in\bbZ_N}\Big(\cos(2\pi nt/N-\varphi)-\frac{a}{2T}\Big)^k<0\qquad\text{($t\in \bbZ_N\ba\{0\}$, $\varphi$ réel),}
\]
donc
\[
P(\text{non}\eqref{enum:4})<N\nu\exp\Big(-\tau N\frac{a^2}{4T^2+2}\Big)\quad(a=\cos\pi/\nu).
\]
\item
Quand $T=1$, la donnée de $\wh x(1)$ détermine $x$. Montrer que cependant on ne peut pas supprimer $\log N$ ou le remplacer par une fonction de $N$ moins vite croissante dans la formule \eqref{eq:4},  \resp \eqref{eq:5} (prendre $N$ pair).
\item
Expliciter la probabilité d'avoir \eqref{enum:1} par le procédé indiqué quand $T=2$.
\end{enumeratea}

\backmatter

\providecommand{\bysame}{\leavevmode ---\ }
\providecommand{\og}{``}
\providecommand{\fg}{''}
\providecommand{\smfandname}{\&}
\providecommand{\smfedsname}{\'eds.}
\providecommand{\smfedname}{\'ed.}
\providecommand{\smfmastersthesisname}{M\'emoire}
\providecommand{\smfphdthesisname}{Th\`ese}

\end{document}